\documentclass[11pt,leqno]{article}
\usepackage{a4wide,amsmath,amssymb,amsfonts,exscale}
\usepackage[curve,matrix,arrow,cmtip]{xy}
\usepackage{hyperref}
\long\def\beginFORGET#1\endFORGET{}

\def\an{{\rm an}}

 \def\mG{{\mathbb G}}

 \def\Aut{\mathop{\rm Aut}\nolimits}

 \def\Gal{\mathop{\rm Gal}\nolimits}
 \def\End{\mathop{\rm End}\nolimits}
 \def\Pic{\mathop{\rm Pic}\nolimits}

 \def\Spec{\mathop{\rm Spec}\nolimits}
 
 \def\deg{\mathop{\rm deg}\nolimits}
 \def\Deg{\mathop{\rm Deg}\nolimits}

\def\Stab{\mathop{\rm Stab}\nolimits}

\def\Supp{\mathop{\rm Supp}\nolimits}

\def\Quot{\mathop{\rm Quot}\nolimits}

\def\diag{\mathop{\rm diag}\nolimits}

\def\GL{\mathop{\rm GL}\nolimits}
\def\SL{\mathop{\rm SL}\nolimits}

\def\ab{{\rm ab}}
\def\ns{{\rm ns}}

\let\phi\varphi
\let\epsilon\varepsilon

\newtheorem{Thm}{Theorem}[section]
\newtheorem{Prop}[Thm]{Proposition}
\newtheorem{Lem}[Thm]{Lemma}

\newtheorem{Cor}[Thm]{Corollary}
\newtheorem{Conj}[Thm]{Conjecture}
\newtheorem{Def}[Thm]{Definition}

\newtheorem{Rem}[Thm]{Remark}

\newenvironment{myequation}
               {\addtocounter{Thm}{1}\begin{equation}}
               {\end{equation}}

\def\qed{{\hskip0pt\unskip\unskip\nobreak\hfil\penalty50
          \hskip1em\hbox{}\nobreak\hfil
           {$\square$}
          \parfillskip=0pt\finalhyphendemerits=0
          \par}\medskip}

\newenvironment{Proof}
               {\noindent{\bf Proof.}\ }
               {\qed}

\newenvironment{Proofof}[1]
               {\noindent{\bf Proof of #1.}\ }
               {\qed}

\newcommand{\BA}{{\mathbb{A}}}

\newcommand{\BC}{{\mathbb{C}}}

\newcommand{\BF}{{\mathbb{F}}}

\newcommand{\BN}{{\mathbb{N}}}

\newcommand{\BP}{{\mathbb{P}}}
\newcommand{\BQ}{{\mathbb{Q}}}

\newcommand{\BZ}{{\mathbb{Z}}}

\newcommand{\Fa}{{\mathfrak{a}}}

\newcommand{\Fc}{{\mathfrak{c}}}

\newcommand{\Fe}{{\mathfrak{e}}}

\newcommand{\Fm}{{\mathfrak{m}}}
\newcommand{\Fn}{{\mathfrak{n}}}

\newcommand{\Fp}{{\mathfrak{p}}}

\newcommand{\FN}{{\mathfrak{N}}}

\newcommand{\FP}{{\mathfrak{P}}}

\newcommand{\CF}{{\cal F}}

\newcommand{\CK}{{\cal K}}
\newcommand{\CL}{{\cal L}}
\newcommand{\CM}{{\cal M}}

\newbox\mybox
\def\arrover#1{\mathrel{
       \setbox\mybox=\hbox spread 1.4em
              {\hfil$\scriptstyle#1$\hfil}
       \vbox{\offinterlineskip\copy\mybox
             \hbox to\wd\mybox{\rightarrowfill}}}}

\def\larrover#1{\mathrel{
       \setbox\mybox=\hbox spread 1.4em
              {\hfil$\scriptstyle#1\vphantom{g}$\hfil}
       \vbox{\offinterlineskip\copy\mybox
             \hbox to\wd\mybox{\leftarrowfill}}}}

\def\ontoover#1{\mathrel{
       \setbox\mybox=\hbox spread 1.4em
              {\hfil$\scriptstyle#1\vphantom{g}$\hfil}
       \vbox{\offinterlineskip\copy\mybox
             \hbox to\wd\mybox{\rightarrowfill\hskip-2.8mm
                               $\rightarrow$}}}}
\def\leftontoover#1{\mathrel{
       \setbox\mybox=\hbox spread 1.4em
              {\hfil$\scriptstyle#1\vphantom{g}$\hfil}
       \vbox{\offinterlineskip\copy\mybox
             \hbox to\wd\mybox{$\leftarrow$\hskip-2.8mm
                               \leftarrowfill}}}}
\let\longto\longrightarrow
\let\into\hookrightarrow
\let\onto\twoheadrightarrow

\def\isoto{\arrover{\sim}}
\def\longinto{\lhook\joinrel\longrightarrow}

\def\invlim{\mathop{\vtop{\hbox{\rm lim}\vskip-8pt
        \hbox{\hskip1pt$\scriptstyle\longleftarrow$}\vskip-1pt}}}

\def\Cinf{{\BC}_\infty}
\def\Kinf{K_\infty}
\def\Ahat{\hat{A}}
\def\Af{\BA_f}
\def\HCM{H_{\rm CM}}

\NoComputerModernTips
\newdir^{ (}{{}*!/-3pt/\dir^{(}}
\newdir_{ (}{{}*!/-3pt/\dir_{(}}

\begin{document}

\title{Special subvarieties of Drinfeld modular varieties}
\author{Florian Breuer (fbreuer@sun.ac.za)}
\maketitle

\begin{abstract}
We explore an analogue of the Andr\'e-Oort conjecture for subvarieties of Drinfeld modular varieties. The conjecture states that a subvariety $X$ of a Drinfeld modular variety contains a Zariski-dense set of complex multiplication (CM) points if and only if $X$ is a ``special'' subvariety (i.e. $X$ is defined by requiring additional endomorphisms). We prove this conjecture in two cases. Firstly when $X$ contains a Zariski-dense set of CM points with a certain behaviour above a fixed prime (which is the case if these CM points lie in one Hecke orbit), and secondly when $X$ is a curve containing infinitely many CM points without any additional assumptions.
\end{abstract}

\noindent{\em 2000 Mathematics Subject Classification: 11G09, 14G35.}

\section{Introduction}

\subsection{The Andr\'e-Oort conjecture}

The Andr\'e-Oort conjecture for complex Shimura varieties states that any irreducible subvariety $X$ of a Shimura variety $S/\BC$ containing a Zariski-dense set of special points must be of Hodge type. For a good overview of this topic, see Rutger Noot's survey paper \cite{Noot04}. 

There has been significant progress on this conjecture recently, which we classify into four flavours. 

{\bf (A)} Products of modular curves. Yves Andr\'e \cite{Andre98} proved the Andr\'e-Oort conjecture when $S$ is a product of two modular curves. Slightly earlier, Bas Edixhoven \cite{Edixhoven98} found a proof assuming the Generalised Riemann Hypothesis (GRH), which generalised well to products of several modular curves \cite{Edixhoven05}, and to products of Shimura curves (by Andrei Yafaev \cite{Yafaev01}). This method, involving Galois orbits of special points and Hecke correspondences, is central to almost all further progress.

{\bf (B)} Very special points. There are several results showing that a subvariety $X\subset S$ is special if $X$ contains a Zariski-dense set of special points which satisfy additional properties. For example, Ben Moonen \cite[\S5]{Moonen98b} treated the case where $S$ is a moduli space of principally polarised abelian varieties, and $X$ contains a Zariski-dense set of special points corresponding to CM abelian varieties which are the canonical lifts of their reductions modulo a fixed prime. This has been generalised to arbitrary Shimura varieties $S$ by Yafaev \cite{YafaevPrep}. 


{\bf (C)} The case when $X$ is a curve. This is the form originally conjectured by Andr\'e in \cite[Problem 9]{Andre89}. When $X$ is a curve on a Hilbert modular surface the conjecture was proved by Edixhoven \cite{Edixhoven01} under GRH. Edixhoven and Yafaev \cite{Edixhoven-Yafaev03} proved the case when $X$ is a curve in a general Shimura variety, but containing infinitely many special points with isomorphic $\BQ$-Hodge structures. The general case for subcurves $X$ was proved by Yafaev \cite{YafaevPrep2}, though again assuming GRH.
%

{\bf (D)} The general case. Recently, Bruno Klingler, Emmanuel Ullmo and Yafaev have announced a proof of the full Andr\'e-Oort conjecture, still assuming GRH, see \cite{KlinglerYafaev, UllmoYafaev}. This builds on equidistribution results of Laurent Clozel and Ullmo \cite{ClozelUllmo}, in addition to Galois actions on subvarieties of Hodge type. 

\subsection{Characteristic $p$}

It is interesting to explore analogues of these developments in characteristic $p$, particularly in the context of Drinfeld modules. Results analogous to (A) above were obtained by the author for products of Drinfeld modular curves in odd characteristic \cite{BreuerPrep, BreuerTrans}. Here special points are called CM points, as they correspond to tuples of rank 2 Drinfeld modules with complex multiplication.

The present paper continues this programme, by addressing analogues of the results in (B) and (C) above. More precisely, suppose $M$ is a Drinfeld modular variety parametrising rank $r$ Drinfeld modules with some level structure. A point on $M$ is called a CM point if the underlying Drinfeld module has complex multiplication, and an irreducible algebraic subvariety $X\subset M$ is called special if it is (roughly) an irreducible component of the locus of Drinfeld modules with additional endomorphisms (see Section \ref{specialsubvarsection} for the full definition). Then we conjecture that an irreducible algebraic subvariety 
$X\subset M$ contains a Zariski dense set of CM points if and only if $X$ is special. 

One difficulty is that much of the machinery needed to study such questions is lacking in the literature, and this article provides some of that machinery. Our techniques allow us to prove two main results (see Section \ref{AOconjSect} for the full statements). First, we prove an analogue of Moonen's result, where the CM points all have good behaviour above a fixed prime. This is the case, in particular, when all the CM points lie in one Hecke orbit. Secondly, we show that an irreducible curve $X\subset M$ contains infinitely many CM points if and only if $X$ is special. In particular, our Conjecture holds for Drinfeld modular surfaces, which parametrise Drinfeld modules of rank three.

Our proof strategy is similar to that developed by Edixhoven and Yafaev, though translation of the details is not automatic. Let $X\subset M$ be an irreducible algebraic subvariety, containing a Zariski-dense set of CM points. We first reduce to the case where $X$ does not lie in a proper special subvariety of $M$ (we say $X$ is Hodge generic). Next we study Galois actions on CM points to show that $X$ is stabilised by a suitable Hecke correspondence (this requires a strong form of the \v{C}ebotarev Density Theorem in case (C), which is where GRH is needed in the classical case). The image under this correspondence is then shown to be irreducible, so that $X$ contains complete Hecke orbits. Finally, one shows that such Hecke orbits are dense in an irreducible component of $M$, and the result follows.   

To tackle the general case (D) will require Galois action on special subvarieties, and equidistribution techniques to handle special subvarieties with bounded Galois orbits. This is the subject of an ongoing PhD thesis, and we will not touch it here.

%
%
%
%
%
%
%
%
%

\subsection{Layout}

Section \ref{Section2} gathers some basic background on Drinfeld modular varieties, including Hecke correspondences and the theory of complex multiplication. Many of these results are already known, and are included, as concisely as possible, to fix our notation.

In Section \ref{Section3} we define special subvarieties of Drinfeld modular varieties, inspired by the definition of subvarieties of Hodge type (Definition \ref{SpecialSubsDef}). 
In \S\ref{AOconjSect} we state an analogue of the Andr\'e-Oort conjecture (Conjecture \ref{AOconj}) and our main results (Theorems \ref{MainMoonen} and \ref{MainCurves}).

Section \ref{Section4} deals with the analytic theory of Drinfeld modular varieties, treating Hecke correspondences in some detail (\S\ref{HCsubsect} and \S\ref{ExplicitSubSect}), and then going on to prove that certain Hecke orbits are Zariski-dense in the underlying moduli space (Theorem \ref{density}).

In Section \ref{section4.5} we exhibit explicit polynomial equations satisfied by Hecke correspondences over rational function fields. We use this to prove an intersection result (Proposition \ref{intersections}) needed in the proof of Theorem \ref{MainCurves}.

In Section \ref{Section5} we show that the images of irreducible subvarieties under certain Hecke correspondences are irreducible (Theorem \ref{irreducible}).

Finally, the main results are proved in Section \ref{MainProofsSect}.

\paragraph{Acknowledgements.} The author is grateful to Bas Edixhoven and Richard Pink for helpful discussions and  detailed readings of earlier versions of this paper.
During the course of this work, the author has also benefited from discussions with Gebhard B\"ockle, Urs Hartl, Marc Hindry, Hans-Georg R\"uck, Andrei Yafaev and Chia-Fu Yu. 

\subsection{Notation and conventions}

All algebraic varieties are assumed to be closed, unless stated otherwise.

The reader is assumed familiar with the basic theory of function field arithmetic and Drinfeld modules.
Standard references for Drinfeld modules and their moduli spaces are \cite{Deligne-Husemoller87}, \cite{Drinfeld74},
\cite{GekelerDMC}, \cite{GossBS} and \cite{Saidi96}.

The following notation will be used throughout the article:
\begin{itemize}
\item $|S|$ denotes the cardinality of the finite set $S$.
\item $K$ is a global function field with exact field of constants $\BF_q$.
\item $\infty$ is a fixed place of $K$.
\item $A$ is the Dedekind domain of elements of $K$ regular away from $\infty$.
\item $|\Fn| := |A/\Fn|$ for a non-zero ideal $\Fn\subset A$.
\item $\deg(a) := \log_q|aA|$ for $a\in A, a\neq 0$.
\item $\Kinf$ is the completion of $K$ at $\infty$.
\item $\Cinf := \hat{\bar{K}}_\infty$ is the completion of an algebraic closure of $\Kinf$.
\item $A_\Fp = \invlim_{n} A/\Fp^n$ is the completion of $A$ at $\Fp$.
\item $K_\Fp = \Quot(A_\Fp)$ is the completion of $K$ at $\Fp$.
\item $\Ahat = \invlim_{\Fn} A/\Fn = \prod_\Fp A_\Fp$ is the profinite completion of $A$. 
\item $\Af = \Ahat\otimes_A K$ is the ring of finite ad\`eles of $K$.
\item $M_A^r(\CK)$ is the moduli scheme of rank $r$ Drinfeld $A$-modules with $\CK$-level structure, 
where $\CK$ is an open subgroup of $\GL_r(\Ahat)$.
\item $M_A^r(\CK)_R := M_A^r(\CK)\times_{\Spec A} \Spec R$ denotes the base extension of $M_A^r(\CK)$ to the $A$-ring $R$.
\end{itemize}

\section{Drinfeld modular varieties}\label{Section2}

\subsection{Drinfeld modules over a scheme}


Let $\gamma_S : S\to\Spec A$ be an $A$-scheme.
Recall (e.g. \cite[Def. 5.1]{Deligne-Husemoller87}) that a Drinfeld $A$-module of rank $r$ over $S$
is a pair $(\CL,\phi)_S$ consisting of a line bundle $\CL/S$ and a ring homomorphism 
$\phi : A\to\End_{\BF_q}(\CL);\; a\mapsto \phi_a$ which over any trivializing open subset $\Spec B\subset S$ is given by
$\phi_a = \sum_{i=0}^m a_i\tau^i\in \End_{\BF_q}(\mG_{{\mathrm a},B}) = B\{\tau\}$, where $\tau : b\mapsto b^q$ is the $q$-power Frobenius, and $a_0 = \gamma_S^\sharp(a)$, $m=r\deg(a)$ and $a_m$ is a unit in $B$. The ring $B\{\tau\}$ is called the ring of twisted polynomials over $B$, and the following commutation relation holds: $\tau b = b^q\tau$ for all $b\in B$.

When $F$ is a field, $S=\Spec(F)$ and $\CL=\mG_{{\mathrm a},F}$, we call $\phi$ a Drinfeld module over $F$. We say it has generic characteristic if $\gamma_S^\sharp : A\hookrightarrow F$.

Let $\Fn\subset A$ be a non-zero ideal, then a level-$\Fn$ structure on $(\CL,\phi)_S$ is a morphism of $A$-modules
$\alpha : (\Fn^{-1}/A)^r \to \CL(S)$ such that $\sum_{x\in (\Fn^{-1}/A)^r}\alpha(x) = \phi[\Fn]$ as Cartier
divisors on $\CL$. A full level structure on $(\CL,\phi)_S$ is a morphism of $A$-modules 
$\alpha:(K/A)^r\to\CL(S)$ such that $\alpha|_{(\Fn^{-1}/A)^r}$ is a level-$\Fn$ structure for every non-zero ideal
$\Fn\subset A$.

A morphism $(\CL,\phi,\alpha)_S \to (\CM,\psi, \beta)_S$ of Drinfeld modules with a given level structure is a morphism $f : \CL \to \CM$ of line bundles over $S$ such that $f\circ \phi_a = \psi_a\circ f$ for all $a\in A$, and $f\circ\alpha = \beta$. It is an isomorphism if it is an isomorphism of line bundles. A non-zero morphism is called an {\em isogeny}.


\subsection{Moduli schemes}\label{IntroModuliSect}

Consider the functor $\CF_A^r(\Fn) : \text{$A$-Schemes} \to \text{Sets}$, which maps $S$ to the set of isomorphism classes of 
rank $r$ Drinfeld $A$-modules over $S$ with level-$\Fn$ structure. When $\Fn$
is divisible by two distinct primes (we call $\Fn$ {\em admissible}), 
then $\CF_A^r(\Fn)$ is representable by a fine moduli scheme $M_A^r(\Fn)$.
Moreover, $M_A^r(\Fn)$ is an affine scheme, flat of dimension $r-1$ over $A$, and is smooth over 
$\Spec A \smallsetminus V(\Fn)$. For any $\Fm\subset\Fn$ one has the forgetful functor $M_A^r(\Fm)\to M_A^r(\Fn)$,
and the projective limit of $M_A^r(\Fn)$ over all admissible $\Fn$ is denoted by $M_A^r$. 
This is the fine moduli scheme parametrizing rank $r$ Drinfeld $A$-modules with full level structures. 

Now $\GL_r(\Af)$ acts from the left on 
$M_A^r$ (see \cite[\S I.3]{GekelerDMC}), and for an open subgroup $\CK\subset\GL_r(\Af)$ we obtain $M_A^r(\CK):=\CK\backslash M_A^r$ as the coarse moduli scheme parametrizing Drinfeld $A$-modules of rank $r$ with level-$\CK$
structures. We see that $M_A^r(\CK(\Fn))=M_A^r(\Fn)$, where $\CK(\Fn):=\ker\big(\GL_r(\Ahat)\to\GL_r(A/\Fn)\big)$, even when $\Fn$ is not admissible. $M_A^r(\CK)$ is a fine moduli scheme
if $\CK\subset\CK(\Fn)$ for some admissible $\Fn\subset A$.
We will write $M_A^r(1) = M_A^r(\GL_r(\Ahat))$ for the coarse moduli scheme parametrizing rank $r$ Drinfeld $A$-modules
without additional level structures.

\subsection{Hecke correspondences}\label{IntroHeckeSect}

Let $g\in\GL_r(\Af)$. The action of $g$ on $M_A^r$ induces a {\em Hecke correspondence} $T_g$ on $M_A^r(\CK)$:
\[
T_g : M_A^r(\CK) \stackrel{\pi}{\longleftarrow} M_A^r \stackrel{g}{\longrightarrow} M_A^r \stackrel{\pi}{\longrightarrow} M_A^r(\CK).
\]
Since the composition $\pi\circ g$ is the quotient by $g^{-1}\CK g$, the Hecke correspondence factors through $M_A^r(\CK_g)$, where $\CK_g := \CK\cap g^{-1}\CK g$,
\[
T_g : M_A^r(\CK) \stackrel{\pi}{\longleftarrow} M_A^r(\CK_g)  \stackrel{\pi\circ g}{\longrightarrow} M_A^r(\CK).
\]

It follows that $T_g$ is a finite algebraic correspondence of degree $\deg(T_g)=\deg(\pi) = [\CK/Z\cap\CK : \CK_g/Z\cap\CK_g]$, where $Z$ denotes the centre of $\GL_r(\Ahat)$. 
Note that for any $x\in M_A^r(\CK)$, the Drinfeld modules corresponding to the points $T_g(x)$ are linked to the Drinfeld module corresponding to $x$ by isogenies specified by the choice of $g$.


Let $\Fn\subset A$ be a non-zero ideal. We denote by $T_\Fn$ the correspondence on $M_A^r(\CK)$ formed by the union of all Hecke correspondences encoding incoming cyclic $\Fn$-isogenies (i.e. if $y\in T_\Fn(x)$ then there exists an underlying isogeny $f : \phi^y \to \phi^x$ with $\ker f \cong A/\Fn$). When $\CK=\GL_r(\Ahat)$, then $T_\Fn = T_g$ where $g=\diag(\hat\Fn,1,\ldots,1)$, where $\hat\Fn\in\Ahat$ denotes a generator of the principal ideal $\Fn\Ahat$. For an open subgroup $\CK\subset\GL_r(\Ahat)$, the graph of $T_\Fn$ in 
$M_A^r(\CK)\times M_A^r(\CK)$ is the preimage of the graph of $T_\Fn$ in $M_A^r(1)\times M_A^r(1)$. We thus see that $T_\Fn$ is a finite union of Hecke correspondences. 

%

\subsection{Complex multiplication}\label{CMsect}

Let $\phi$ be a rank $r$ Drinfeld $A$-module over a field in generic characteristic. Then the endomorphism ring $R=\End(\phi)$ is a commutative 
$A$-algebra which has rank dividing $r$ as a projective $A$-module. 
Moreover, $R$ is an order in its quotient field $K'$, and the extension $K'/K$ is 
{\em purely imaginary}, which means that $K'$ has only one prime lying above $\infty$, which we again denote by $\infty$.
When $[K':K]=[R:A]=r$, then we say $\phi$ has {\em complex multiplication} (CM) by $R$. 

%

Suppose now that $\phi$ has CM by the ring $R$. Then we may view $\phi$ as a rank one Drinfeld $R$-module, even if $R$ is not integrally closed, thanks to the work of Hayes \cite{Hayes79}. Denote by $H_R$ the {\em ring class field} associated to the order $R$
(in general, if $F$ is a finite extension of $K$, we define the {\em Hilbert class field of $F$} to be the ring class field of the integral closure of $A$ in $F$). Then $H_R/K'$ is an abelian extension, totally split at $\infty$ and unrammified outside the conductor of $R$ in $A'$, and $\Gal(H_R/K')\cong\Pic(R)$. 

Hayes showed in \cite[\S8]{Hayes79} that up to isomorphism, there are  precisely $|\Pic(R)|$ rank one Drinfeld $R$-modules, each defined over the field $H_R$, and $\Gal(H_R/K')$ acts on these isomorphism classes like isogenies. More precisely, let $\Fa\subset R$ be a non-zero ideal, let $\sigma_{\Fa}=(\Fa,H_R/K')\in\Gal(H_R/K')\cong\Pic(R)$ correspond to the class of $\Fa$, and let $\phi$ be a rank one Drinfeld $R$-module. Then $\phi^{\sigma_{\Fa}} = \Fa * \phi$. In other words, there is an isogeny $\phi \rightarrow \phi^{\sigma_{\Fa}}$ with kernel $R/\Fa$.

A point $x\in M_A^r(\CK)(\Cinf)$ is called a {\em CM point} if the corresponding Drinfeld module $\phi^x$ (determined by $x$ up to isomorphism) has complex multiplication.

\begin{Prop}\label{CMmain}
Let $\CK\subset\GL_r(\Ahat)$ be an open subgroup, and $x\in M_A^r(\CK)(\Cinf)$ a CM point with corresponding endomorphism ring $R=\End(\phi^x)$, which is an order in the field $K'=\Quot(R)$. 
\begin{enumerate}
\item The field of definition $K'(x)$ of $x$ over $K'$ is a finite abelian extension of $K'$ containing the ring class field $H_R$ of $R$.
\item Let $\FN\subset R$ be an ideal such that $R/\FN \cong A/\Fn$, where $\Fn = \FN\cap A$. Denote by $\sigma = (\FN, K'(x)/K')^{-1}$ the inverse Artin element. Then $\sigma(x)\in T_\Fn(x)$.
%
%
\end{enumerate}
\end{Prop}

\begin{Proof}
{\bf (1)}
We may choose $\phi^x$ defined over $H_R$. The full torsion submodule (and hence a full level structure) of $\phi^x$ is then defined over ${K'}^\ab$, the maximal abelian extension of $K'$ which splits completely at $\infty$. It follows that any point of $M_A^r(\Cinf)$ lying above $x$ is defined over ${K'}^\ab$, hence so is $x$ itself. The $\Gal({K'}^\ab/H_R)$-orbit of $x$ lies completely above the point of $M_A^r(1)(\Cinf)$ corresponding to $\phi^x$, and it follows that $[K'(x):H_R]$ can be no larger than the degree of the morphism $M_A^r(\CK)\to M_A^r(1)$. 

{\bf (2)}
There is an isogeny of rank 1 Drinfeld $R$-modules $\FN^{-1} * \phi^x = (\phi^x)^\sigma\to\phi^x$ with kernel isomorphic to $R/\FN$, which is automatically also an isogeny of rank $r$ Drinfeld $A$-modules with kernel isomorphic to $A/\Fn$. Note that $A\into R \onto R/\FN$ is a morphism of $A$-modules with kernel $\Fn$, so $R/\FN\cong A/\Fn$ as rings implies the isomorphism for $A$-modules, too.
\end{Proof}

\subsection{CM heights}\label{CMheightsSect}

Let $\phi$ be a rank $r$ Drinfeld $A$-module with complex multiplication by an order $R$ in $K'$, as above. 
Let $\Fc = \{a\in A' \;|\; aA'\subset R\}$ be the conductor of $R$ in $A'$, it has the property that 
$|A'/\Fc|/|R/\Fc| = |A'/R|$. Write $|\Fc|=|A'/\Fc|$, and denote the genus of $K'$ by $g(K')$. 
Then we define the {\em CM height} of $\phi$ by
\begin{myequation}
\HCM(\phi) := q^{g(K')}\cdot|\Fc|^{1/r}.
\end{myequation}
Similarly, for any CM point $x\in M_A^r(\CK)(\Cinf)$ with underlying Drinfeld module $\phi^x$, we set $\HCM(x):=\HCM(\phi^x)$.

\begin{Prop}\label{CMheights}
Let $\CK\subset\GL_r(\Ahat)$ be an open subgroup.
\begin{enumerate}
\item Let $B>0$ be any given constant, then $M_A^r(\CK)(\Cinf)$ contains only finitely many CM points of CM height bounded by $B$.

\item For every $\epsilon>0$ there is a computable constant $C_\epsilon>0$ such that the following holds.
Let $\phi$ be a Drinfeld module with complex multiplication by $R$ as above, then  
\[
|\Pic(R)| > C_\epsilon\HCM(\phi)^{1-\epsilon}.
\]
\end{enumerate}
\end{Prop}

\begin{Proof}
{\bf (1)} For a given $g\geq 0$ the field $K$ has only finitely many extensions $K'$ of degree $r$ with constant field contained in $\BF_{q^r}$ and genus $g(K')\leq g$, by \cite[Theorem 8.23.5]{GossBS}. 
For each such extension, let $A'$ be the integral closure of $A$ in $K'$. For each ideal $\Fc\subset A'$, there are only finitely many orders $R\subset A'$ of conductor $\Fc$ which are projective $A$-modules of rank $r$, since $A'$ has only finitely many $A$-submodules of index $\leq |A'/\Fc|$. 
For any such order $R$, there are at most $|\Pic(R)|$ rank $r$ Drinfeld $A$-modules with complex multiplication by $R$, and each such Drinfeld module corresponds to only finitely many CM points on $M_A^r(\CK)(\Cinf)$.

{\bf (2)} 
We start with the following lower bound for the class number of $K'$ \cite[Prop. 3.1]{BreuerPrep}.
\[
|\Pic(A')| \geq h(K') \geq \frac{(q-1)(q^{2g(K')}-2g(K')q^{g(K')}+1)}{2g(K')(q^{g(K')+1}-1)}.
\]
Next, we have the exact sequence \cite[\S I.12]{Neukirch}
\[
1\to A'^\times/R^\times \longto (A'/\Fc)^\times/(R/\Fc)^\times \longto \Pic(R) \longto \Pic(A') \to 1
\]
and since $|A'^\times/R^\times|\leq q^r$ it remains to bound $(A'/\Fc)^\times/(R/\Fc)^\times$ from below in terms of $|\Fc|$. 

By the Invariant Factor Theorem \cite[\S22]{CurtisReiner} there exist $A$-ideals $\Fa_1,\ldots,\Fa_r$ and $\Fe_1|\Fe_2|\cdots|\Fe_r$, and elements $m_1,\ldots,m_r\in A'$ such that, as $A$-modules,
\[
A'=\Fa_1m_1\oplus\cdots\oplus\Fa_rm_r, \quad\text{and}\quad R=\Fa_1\Fe_1m_1\oplus\cdots\oplus\Fa_r\Fe_rm_r.
\]
The conductor of $R$ satisfies $\Fc\cap A = \Fe_r$, and we have
\[
\frac{|A'/\Fc|}{|R/\Fc|} = |A'/R| =
\left|\frac{\Fa_1m_1\oplus\cdots\oplus\Fa_rm_r}{\Fa_1\Fe_1m_1\oplus\cdots\oplus\Fa_r\Fe_rm_r}\right|
\geq |A/\Fe_r| \geq |A'/\Fc|^{1/r}.
\]

On the other hand

\[
\frac{|(A'/\Fc)^\times|}{|(R/\Fc)^\times|} \geq \frac{|(A'/\Fc)^\times|}{|R/\Fc|} =
\frac{|A'/\Fc|}{|R/\Fc|}\prod_{\Fp|\Fc}\left(1-\frac{1}{|\Fp|}\right)
\geq C_\varepsilon|\Fc|^{1/r-\varepsilon}.
\]
%
%
%
%
The result follows.
\end{Proof}

\section{Special subvarieties}\label{Section3}

\subsection{Enlarging the base ring}

Let $r'|r$ be an integer, and let $K'/K$ be a purely imaginary extension of degree $r/r'$.
Let $A'$ denote the integral closure of $A$ in $K'$. Let $\Fn\subset A$ be a non-zero ideal, and set $\Fn'=\Fn A'$. Note that $A'/\Fn' \cong (A/\Fn)^{r/r'}$ as $A$-modules.
We fix an isomorphism of $A$-modules $(\Fn^{-1}/A)^r \isoto (\Fn'^{-1}/A')^{r'}$.

Now let $S\to\Spec(A')$ be an $A'$-scheme and
let $(\CL,\phi',\alpha')_S$ be a rank $r'$ Drinfeld $A'$-module with level-$\Fn'$ structure. Then restricting
$\phi'$ to $A$, we obtain a rank $r$ Drinfeld $A$-module $(\CL,\phi,\alpha)_S$, where $\phi=\phi'|_A$ and 
\[
\alpha : (\Fn^{-1}/A)^r \isoto (\Fn'^{-1}/A')^{r'} \arrover{\alpha'} \CL(S)
\]
is a level-$\Fn$ structure on $\phi$.

We thus obtain a transformation of functors $F_{A'/A}$ from $\CF_{A'}^{r'}(\Fn')$ to $\CF_A^r(\Fn)$ (restricted to the subcategory of $A'$-schemes), and correspondingly a 
canonical morphism of moduli schemes $f_{A'/A} : M_{A'}^{r'}(\Fn')\to M_A^r(\Fn)_{A'}$.  

\begin{Lem}\label{injectivity}
Fix embeddings $A \subset A' \subset \Cinf$, then the above canonical morphism is injective on $\Cinf$-valued points:
\[ 
f_{A'/A} : M_{A'}^{r'}(\Fn')(\Cinf)\longinto M_A^r(\Fn)(\Cinf) 
\]
\end{Lem}

\begin{Proof}
Two Drinfeld $A'$-modules over $\Cinf$ correspond to $A'$-lattices $\Lambda_1,\Lambda_2\subset\Cinf$. They are isomorphic as Drinfeld $A$-modules if and only if $\Lambda_1=c\Lambda_2$ for some $c\in\Cinf$, in which case they are also isomorphic as Drinfeld $A'$-modules.
%
\end{Proof}

\subsection{Properness}

\begin{Prop}\label{properness}
Let $r',K',A',\Fn$ and $\Fn'$ be as above.
Suppose that $\Fn$ is admissible, so that $\Fn'$ is admissible and  $M_A^r(\Fn)$ and $M_{A'}^{r'}(\Fn')$ are fine moduli schemes.
Then the morphism $f_{A'/A} : M_{A'}^{r'}(\Fn')_{K'}\to M_A^r(\Fn)_{K'}$ is proper.
\end{Prop}

\begin{Proof}
The morphism is separated and of finite type, as both schemes are spectra of $K'$-algebras of finite type. 

We use the valuative criterion for properness. Let $R$ be a discrete valuation ring with quotient field $L$ 
and valuation $v$, provided with a map $A'\into R$ (base extension to $K'$ forces generic characteristic). Given a commutative diagram
\[
\xymatrix{
\Spec L \ar[r]^f\ar[d] & M_{A'}^{r'}(\Fn')_{K'}\ar[d] \\
\Spec R \ar[r]^g\ar@{-->}[ur]^{\exists!h} & M_A^r(\Fn)_{K'}
}
\]
with given morphisms $f$ and $g$, we want to show that there exists a morphism $h$ (we already have uniqueness).
The morphism $f$ corresponds to an isomorphism class of rank $r'$ Drinfeld $A'$-modules $\phi : A'\to L\{\tau\}$, and
$g$ corresponds to an isomorphism class of rank $r$ Drinfeld $A$-modules $\psi : A \to R\{\tau\}$, such that the leading coefficient of any $\psi_a$ is a unit in $R$ (so $\psi$ is a Drinfeld module over $\Spec R$). That the diagram commutes means that $\phi|_A \cong \psi$. By choosing suitable representatives in the isomorphism classes, we may assume that
$\phi|_A=\psi$.
To construct the morphism $h$, we will show that $\phi$ is already defined 
over $R$, and that the leading coefficient of each $\phi_a$, $a\in A'$, is a unit.

Let $0\neq a\in A'$ with minimal polynomial 
$f(X)=b_0 + b_1X + \cdots + b_{n-1}X^{n-1} + X^n\in A[X]$.
Then
\begin{myequation}\label{eq1}
0=f^\phi(\phi_a) = \phi_{b_0} + \phi_{b_1}\phi_a + \cdots + \phi_{b_{n-1}}\phi_a^{n-1} + \phi_a^n.
\end{myequation} 

Let $0\neq\lambda\in\bar{L}$ such that $\phi_a(\lambda)=0$. Applying the right hand side of (\ref{eq1}) to $\lambda$, 
we get $\phi_{b_0}(\lambda) = 0$. As $b_0\in A$, $\ker\phi_{b_0}$ is \'etale and $v(\lambda)=0$. 
As $\phi_a$ is separable we can write 
\[
\phi_a(T) = a_0T\cdot\prod_{0\neq\lambda\in\ker\phi_a}\left(1-\frac{T}{\lambda}\right)
\quad\text{with $a_0\neq 0$, and $v(\lambda)=0$ for all $\lambda$.}
\]
It follows that $\phi_a(T)\in R[T]$. It remains for us to show that the leading coefficient of $\phi_a$ is a unit.

Write $\phi_a = a_0 + a_1\tau + \cdots a_m\tau^m$. Multiplying out the right hand side of (\ref{eq1}), we find that
the coefficient of the highest power of $\tau$ is of the form
\[
\sum_{i\in I}b'_i a_m^{\left(\frac{q^{mi}-1}{q^m-1}\right)q^{d_i}} = 0
\]
where $I\subset\{0,\ldots,n\}$, and each $b_i'\in R^\times$ is the leading coefficient of $\phi_{b_i}$ ($b_n:=1$), which has
degree $d_i=r\deg(b_i)$. Now, we have at least two terms, and the exponent of $a_m$ is different
in every term (the $i$th exponent has $i-1$ non-zero digits in base $q$). Therefore we must have $v(a_m)=0$, which is
what we wanted to prove.

Lastly, since $\phi[\Fn]=\psi[\Fn']$ is \'etale over $\Spec R$ we see that the level structures are automatically compatible, and thus we have a well-defined morphism $h:\Spec R\to M_{A'}^{r'}(\Fn')_{K'}$.
%
\end{Proof}

This result has been generalized to Abelian sheaves and Drinfeld-Anderson Shtuka by Urs Hartl and Markus Hendler \cite{Hartl}.

\subsection{Special subvarieties}
\label{specialsubvarsection}

%

Proposition \ref{properness} shows that $f_{A'/A}\big(M_{A'}^{r'}(\Fn')_{K'}\big)$ is a closed subvariety of $M_A^r(\Fn)_{K'}$. By a mild (by Proposition \ref{injectivity}) abuse of notation, we denote this image again by $M_{A'}^{r'}(\Fn')_{K'}$. It is
the locus of those Drinfeld modules with endomorphism ring containing $A'$ (via the fixed embedding $A\into A'$).

We now define the special subvarieties of a Drinfeld modular variety.

\begin{Def}\label{SpecialSubsDef}
{\em
Let $\Fn\subset A$ be an admissible ideal.
Let $X\subset M_A^r(\Fn)_{\Cinf}$ be an irreducible subvariety. 
Then $X$ is a {\em special} subvariety if $X$ is an irreducible component of $T_g(M_{A'}^{r'}(\Fn')_{\Cinf})$ for some $g\in\GL_r(\Af)$, $r'|r$ and $A'/A$.

An irreducible subvariety $X\subset M_A^r(1)_{\Cinf}$ is {\em special} if it is the image of a special subvariety under the canonical projection $M_A^r(\Fn)_{\Cinf} \to M_A^r(1)_{\Cinf}$, for some admissible $\Fn\subset A$.

Let $\CK\subset\GL_r(\Ahat)$. An irreducible subvariety $X\subset M_A^r(\CK)_{\Cinf}$ is {\em special} if its image under the canonical projection $M_A^r(\CK)_{\Cinf} \to M_A^r(1)_{\Cinf}$ is special.
}
\end{Def}

Specialness is thus a property of endomorphism rings only, and does not depend on level structures. We now show that in a sense Special subvarieties are the irreducible components of maximal subvarieties parametrising Drinfeld modules with prescribed endomorphisms. 

Let $R$ be a ring that occurs as endomorphism ring for some rank $r$ Drinfeld $A$-module in generic characteristic, it is an order of conductor $\Fc$ in an imaginary extension $K'/K$ of degree $r'|r$. Denote by $A'$ the integral closure of $A$ in $K'$. Let $\Fn\subset A$ be an admissible ideal, and suppose that $\Fc$ divides $\Fn':=\Fn A'$. Let $X\subset M_A^r(\Fn)_{\Cinf}$ be an irreducible algebraic subvariety, and denote by $\phi^X$ the Drinfeld module over $X$ defined by $X\into M_A^r(\Fn)_{\Cinf}$, it has level structure $\alpha_X : (\Fn^{-1}/A)^r \stackrel{\sim}{\longrightarrow} \Phi^X[\Fn]$. Suppose we have a fixed embedding $R\into \End(\phi^X)$.

Denote by $\phi^X[\Fc]$ the group scheme of $\Fc$-torsion points of $\phi^X$, and $H:=\alpha_X^{-1}(\phi^X[\Fc])\subset (\Fn^{-1}/A)^r$. For a geometric point $x\in X(\Cinf)$ we denote by $\phi^x$ the pull-back of $\phi^X$ to $x$ with level structure $\alpha_x$, then we have $R\into\End(\phi^x)$ and $\alpha_x(H)=\phi^x[\Fc]$. Let $h\in\GL_r(\Af)$ be an element whose kernel on $(\Fn^{-1}/A)^r$ is $H$, so $x\mapsto y:=h(x)$ encodes the isogeny $\phi^x \to \phi^y \cong \phi^x/\phi^x[\Fc]$. By \cite[Prop. 4.7.19]{GossBS} we have $A'\into\End(\phi^y)$. Moreover, since $\phi^x[\Fc]\subset\phi^x[\Fn]$ the group $\CK(\Fn)$ acts trivially on $H$, and so every branch of the Hecke correspondence $T_h$ encodes this same underlying isogeny. In other words, $A'\into\End(\phi^y)$ for every $y\in T_h(x)$.

Thus for any irreducible component $Y\subset T_h(X)$ we have $A'\into \End(\phi^Y)$, and consequently $Y\subset f_{A'/A}(M_{A'}^{r'}(\Fn')_{\Cinf}))$ for the given embedding $A\into A'$. $X$ is special if $Y=f_{A'/A}(M_{A'}^{r'}(\Fn')_{\Cinf}))$, so we see that $X$ is special if $X$ is maximal with respect to the conditions: (1) $X\subset M_A^r(\Fn)_{\Cinf}$ is closed and irreducible, and (2) $R\into\End(\phi^X)$. Conversely, special subvarieties are maximal in this respect. 
 In particular, CM points are the special subvarieties of dimension zero.

Now let $X\subset M_A^r(\CK)_{\Cinf}$ be an irreducible algebraic subvariety, where $\CK\subset\GL_r(\Ahat)$ is open. Denote by $X'\subset M_A^r(\Fn)_{\Cinf}$ an irreducible preimage of the image of $X$ in $M_A^r(1)_{\Cinf}$, for $\Fn\subset A$ as above. Let $A'/A$ be a maximal extension for which the above holds, i.e. $Y\subset T_h(X')$ with $\End(\phi^Y)=A'$. We call $(A',h)$ a {\em Hodge pair for $X$}. If $A'=A$, then we say $X$ is {\em Hodge generic}. $X$ is Hodge generic if and only if it contains a point $x\in X(\Cinf)$ with $\End(\phi^x)=A$.


We remark that special subvarieties of Drinfeld modular varieties are analogous to subvarieties of PEL type in moduli spaces of abelian varieties over $\BC$, even though our definition mimicks that of subvarieties of Hodge type.
However, moduli spaces of complex abelian varieties contain subvarieties of Hodge type which are not of PEL type, the situation thus being more complicated than in the case of Drinfeld modular varieties.

\subsection{The Andr\'e-Oort Conjecture for Drinfeld modular varieties}\label{AOconjSect}

We are now ready to state our analogue of the Andr\'e-Oort Conjecture for Drinfeld modular varieties.

\begin{Conj}\label{AOconj}
Let $X\subset M_A^r(\CK)_{\Cinf}$ be an irreducible subvariety. 
Then $X(\Cinf)$ contains a Zariski-dense set of CM points if and only if $X$ is special.
\end{Conj}

Notice that if $r$ is prime, then the only special subvarieties of $M_A^r(\CK)_{\Cinf}$ are CM points and the irreducible components of $M_A^r(\CK)_{\Cinf}$ itself. Thus in this case any infinite set of CM points is expected to be Zariski-dense in an irreducible component of $M_A^r(\CK)_{\Cinf}$. This may be compared to the distribution of torsion points on a simple abelian variety over $\BC$.

It is relatively easy to show that any special subvariety of a Drinfeld modular variety contains a Zariski-dense set of CM points (see Corollary \ref{CMdensity} below), it is the converse that seems to be harder.

A result similar to Conjecture \ref{AOconj} is known for $M_A^r(\CK)$ replaced by the product of Drinfeld modular curves in odd characteristic (see \cite{BreuerPrep} for the case $A=\BF_q[T]$ and \cite{BreuerTrans} for general~$A$).

\begin{Def}\label{def2.1}
{\em
Let $M/L$ be an algebraic field extension. A prime $\Fp$ of $L$ is called {\em residual} in $M$ if there exists a prime $\FP$ of $M$ above $\Fp$ with residual degree $f(\FP|\Fp)=1$. If $R$ is an order in $M$, then $\Fp$ is residual in $R$ if it is residual in $M$ and $\Fp$ does not divide the conductor of $R$. 
}
\end{Def}

The main results of this paper are the following.


\begin{Thm}\label{MainMoonen}
Let $X\subset M_A^r(\CK)_{\Cinf}$ be an irreducible algebraic subvariety. Suppose that $X(\Cinf)$ contains a Zariski-dense subset $\Sigma$ of CM points satisfying the following conditions:
\begin{enumerate}
\item There is a non-zero prime $\Fp\subset A$ such that, for all $x\in\Sigma$, $\Fp$ is residual in the quotient field of $\End(\phi^x)$, \quad and
\item There is an integer $m\in\BN$ such that, for all $x\in\Sigma$, $\Fp^m$ does not divide the conductor of $\End(\phi^x)$.
\end{enumerate}
Then $X$ is special. In particular, if $\Sigma$ lies in one Hecke orbit, then $X$ is special.
\end{Thm}



\begin{Thm}\label{MainCurves}
Let $X\subset M_A^r(\CK)_{\Cinf}$ be an irreducible subcurve. 
Then $X$ contains infinitely many CM points if and only if $X$ is special. 
In particular, Conjecture \ref{AOconj} holds for $r=3$.
\end{Thm}

Theorem \ref{MainMoonen} is an analogue of a result of Ben Moonen \cite[\S5]{Moonen98b} (see also \cite{YafaevPrep}), 
whereas Theorem \ref{MainCurves} is an analogue of Yves Andr\'e's original conjecture on special points on subcurves of Shimura varieties \cite[Problem 9]{Andre89}, proved under GRH by Yafaev in \cite{YafaevPrep2}.
 
The methods used in this paper are an adaptation of the method pioneered by Edixhoven and Yafaev in \cite{Edixhoven98, Edixhoven01, Edixhoven-Yafaev03, YafaevPrep}. The proof of Theorem \ref{MainMoonen} requires us to pay particular attention to level structures, which unfortunately makes the exposition more technical than it might have been. Theorem \ref{MainCurves} requires some intersection theory (Section \ref{section4.5}) which provided the motivation for the paper \cite{BR}. It also requires a strong form of the \v{C}ebotarev Density Theorem, which is where GRH is needed in the classical case. Our result is unconditional, because the GRH holds over function fields.

\section{Analytic theory}\label{Section4}

\subsection{Moduli spaces}\label{modulisection}

Drinfeld $A$-modules over $\Cinf$ correspond to $A$-lattices in $\Cinf$, and the ad\`elic description of such lattices leads to the following analytic description of Drinfeld modular varieties, see e.g. \cite{Deligne-Husemoller87}.
We let
\[
\Omega^r := \BP^{r-1}(\Cinf)\smallsetminus\{\text{Linear subspaces defined over $\Kinf$}\}
\]
denote {\em Drinfeld's upper half-space}, on which $\GL_r(\Kinf)$ acts from the left. Then we have an isomorphism of rigid analytic spaces
%
\[
M_A^r(\CK)^\an(\Cinf) \cong \GL_r(K)\backslash\Omega^r\times\GL_r(\Af)/\CK.
\]
%
Here, $\GL_r(K)$ acts on $\Omega^r$ and on $\GL_r(\Af)$ from the left as usual, and $\CK$ acts on $\GL_r(\Af)$ from the right (and trivially on $\Omega^r$). 

Choose a set of representatives $S$ in $\GL_r(\Af)$ for the finite set $\GL_r(K)\backslash\GL_r(\Af)/\CK$. 
Let $[\omega,h]\in\GL_r(K)\backslash\Omega^r\times\GL_r(\Af)/\CK$, and let $s_h\in S$ be the representative of $h$ in
$\GL_r(K)\backslash\GL_r(\Af)/\CK$. 
Then there exists $(f_h,k_h)\in\GL_r(K)\times\CK$ such that $f_hhk_h=s_h$ and $[\omega,h]=[f_h(\omega),s_h]$. 
It follows that we have a bijection (in fact, a rigid analytic isomorphism)
\begin{eqnarray*}
\GL_r(K)\backslash\Omega^r\times\GL_r(\Af)/\CK & \arrover{\sim} & \coprod_{s\in S}\Gamma_s\backslash\Omega^r \\
\;[\omega,h] & \longmapsto & [f_h(\omega)]_{s_h},\nonumber
\end{eqnarray*}
where for each $s\in S$ we define the arithmetic group $\Gamma_s := s\CK s^{-1}\cap\GL_r(K)$.
We see that the irreducible components of $M_A^r(\CK)^\an(\Cinf)$ (and thus also of $M_A^r(\CK)_{\Cinf}$, by 
\cite[Korollar 3.5]{Lut}) are in bijection with $S$.


The determinant gives a bijection (see for example \cite[II.1.4]{GekelerDMC})
\[
\GL_r(K)\backslash\GL_r(\Af)/\CK \arrover{\sim} K^\times\backslash\Af^\times/\det(\CK).
\]
%
The group on the right hand side is the class group of $A$ corresponding to $\det(\CK)$, and we also notice that 
\[
 K^\times\backslash\Af^\times/\det(\CK)\cong  M_A^1(\det(\CK))^\an(\Cinf),
\]
since $\Omega^1$ is a single point. 
The map sending $[\omega,h]$ to $[\det(h)]$ defines a morphism
\[
M_A^r(\CK)^\an(\Cinf) \longto M_A^1(\det(\CK))^\an(\Cinf)
\]
whose fibres are the irreducible components of $M_A^r(\CK)_{\Cinf}$.

\subsection{Hecke correspondences}\label{HCsubsect}

Let $g\in\GL_r(\Af)$,
then $g$ acts from the left on $\Omega^r\times\GL_r(\Af)$ via $g\cdot(\omega,h) := (\omega,hg^{-1})$. 
Denote by $\pi:\Omega^r\times\GL_r(\Af) \longto \GL_r(K)\backslash\Omega^r\times\GL_r(\Af)/\CK$ the canonical projection.
Then the Hecke correspondence $T_g$ (cf. \S\ref{IntroHeckeSect}) on $M_A^r(\CK)(\Cinf)$ is given by $T_g = \pi\circ g \circ\pi^{-1}$:
\[
\xymatrix{
\Omega^r\times\GL_r(\Af)\ar[r]^{g\cdot}\ar[d]^\pi & \Omega^r\times\GL_r(\Af)\ar[d]^\pi \\
\GL_r(K)\backslash\Omega^r\times\GL_r(\Af)/\CK & \GL_r(K)\backslash\Omega^r\times\GL_r(\Af)/\CK.
}
\]
One easily verifies that this correspondence factors through $\GL_r(K)\backslash\Omega^r\times\GL_r(\Af)/\CK_g$,
where we recall $\CK_g=\CK\cap g^{-1}\CK g$.

If $\det(\CK_g)=\det(\CK)$ then the canonical projection $M_A^r(\CK_g)\to M_A^r(\CK)$ induces a bijection between the sets of irreducible components of $M_A^r(\CK_g)_{\Cinf}$ and $M_A^r(\CK)_{\Cinf}$, respectively, and in this case we say that $T_g$ is {\em irreducible on $M_A^r(\CK_g)_{\Cinf}$}.

In general, the restriction of $T_g$ on any irreducible component of $M_A^r(\CK)_{\Cinf}$ has $\big(\det(\CK):\det(\CK_g)\big)$ irreducible components as an algebraic correspondence.

%

Let $\Fn\subset A$ be a non-zero ideal, and $\hat\Fn\in\Ahat$ a generator of the principal ideal $\Fn\Ahat$. 
Let $\CK=\GL_r(\Ahat)$ and $g=\diag(\hat\Fn,1,\ldots,1)\in\GL_r(\Af)$. Then
\[
\CK_g=\CK_0(\Fn):= \{(a_{ij})\in\GL_r(\Ahat) \;|\; a_{2,1},a_{3,1},\ldots,a_{r,1}\in\Fn\Ahat\},
\]
and clearly $\det(\CK_0(\Fn))=\det(\GL_r(\Ahat))=\Ahat^\times$, so $T_\Fn=T_g$ is irreducible on $M_A^r(1)$.

\medskip

We point out that $M_A^r(\CK_g)_{\Cinf}$, and thus $T_g$, are actually defined over $K$ for all $g\in\GL_r(\Af)$. (Their $\Cinf$-irreducible components are defined over abelian extensions of $K$, a fact we will not need here).

\subsection{Explicit action of $T_g$}\label{ExplicitSubSect}

We now describe the action of $T_g$ on $\coprod_{s\in S}\Gamma_s\backslash\Omega^r$ explicitly. 
Let $[\omega]_s\in\Gamma_s\backslash\Omega^r$, this is the same as 
$[\omega,s]\in \GL_r(K)\backslash\Omega^r\times\GL_r(\Af)/\CK$. It lifts to the set 
\[
\big\{[\omega,sk]\in \GL_r(K)\backslash\Omega^r\times\GL_r(\Af)/\CK_g \;|\; \text{for representatives $k\in\CK$ of $\CK/\CK_g$}\big\}
\]
which in turn maps to the image of $T_g$:
\[
\big\{[\omega,skg^{-1}]\in \GL_r(K)\backslash\Omega^r\times\GL_r(\Af)/\CK \;|\; \text{for representatives $k\in\CK$ of $\CK/\CK_g$}\big\}.
\]
These points all lie in the same irreducible component of $M_A^r(\CK)_{\Cinf}$ since $[\det(sk_1g^{-1})]=[\det(sk_2g^{-1})]$ in 
$K^\times\backslash\Af^\times/\det(\CK)$. Let $s_{g,s}\in S$ correspond to this component, i.e. 
\[
s_{g,s}=f_{s,k,g}\cdot skg^{-1} \cdot k' \in S
\]
for suitable $f_{s,k,g}\in\GL_r(K)$ and $k'\in\CK$.

It follows that $T_g$ maps the $s$-component to the $s_{g,s}$-component in $M_A^r(\CK)_{\Cinf}$, and is given explicitly by
\[
T_g : [\omega]_s \longmapsto \big\{ [f_{s,k,g}(\omega)]_{s_{g,s}} \;|\; \text{for representatives $k\in\CK$ of $\CK/\CK_g$}\big\}.
\]

Now, let $[s_1],\ldots,[s_n]\in\GL_r(K)\backslash\GL_r(\Af)/\CK_g$ lie above $[s]\in\GL_r(K)\backslash\GL_r(\Af)/\CK$, where $n=\big(\det(\CK) : \det(\CK_g)\big)$ is the number of irreducible components of $T_g$. Set $\Gamma_{s_i,g} := s_i\CK_g s_i^{-1}\cap\GL_r(K)$ for $i=1,\ldots,n$, then the components $\Gamma_{s_i,g}\backslash\Omega^r$ of $M_A^r(\CK_g)_{\Cinf}$ lie above the component $\Gamma_s\backslash\Omega^r$ of $M_A^r(\CK)_{\Cinf}$. The $i$th component of $T_g$ then acts explicitly by
\begin{myequation}\label{HeckeDiagram}
\xymatrix{
\Gamma_{s_i,g}\backslash\Omega^r\ar[d]_{\pi}\ar[dr]^{\pi\circ f_{s,k,g}} & \\
\Gamma_s\backslash\Omega^r & \Gamma_{s_{g,s}}\backslash\Omega^r,
}
\end{myequation}
%
\[
T_g : [\omega]_s \longmapsto \big\{ [f_{s,k,g}(\omega)]_{s_{g,s}} \;|\; \text{$[k]\in\CK/\CK_g$ satisfying $[sk]=[s_i]$ in $\GL_r(K)\backslash\GL_r(\Af)/\CK_g$}\big\}.
\]
\subsection{Density of Hecke orbits}

A {\em general Hecke correspondence} $T$ on $M_A^r(\CK)_{\Cinf}$ is a finite algebraic correspondence (not necessarily irreducible), whose action on points is given explicitly by 
\[
[\omega]_{s_1} \mapsto \bigcup_t \{[t(\omega)]_{s_2}\},
\]
where $t\in\GL_r(K)$ ranges over a finite set of matrices, depending on $s_1,s_2\in S$. 

Its inverse is defined by 
$y\in T^{-1}(x)\Leftrightarrow x\in T(y)$. Its $n$th iterate ($n\in\BZ$, iterations of $T^{-1}$ if $n<0$) is denoted $T^n$.
Its orbit of a point $x\in M_A^r(\CK)(\Cinf)$ is denoted $T^\infty(x) := \cup_{n\in\BZ}T^n(x)$. 

%

\begin{Thm}\label{density}
Let $\CK\subset\GL_r(\Ahat)$ be an open subgroup and let $\Fn = \langle N \rangle \subset A$ be a principal ideal, with $|N|>1$.

Let $x\in M_A^r(\CK)(\Cinf)$.
Let $T$ be a general Hecke correspondence preserving the irreducible component of $M_A^r(\CK)_{\Cinf}$ containing $x$ and encoding at least one cyclic $\Fn$-isogeny on this component. 
Then the orbit $T^\infty(x)$ is Zariski-dense in this component.
\end{Thm}

\begin{Proof}
Denote by $\Gamma\backslash\Omega^r$ an irreducible component of $M_A^r(\CK)_{\Cinf}$ containing a point $x=[\omega]$.
The Drinfeld module $\phi^x$ corresponds to a lattice
\[
\Lambda_\omega = \omega_1 A \oplus \cdots \oplus \omega_{r-1} A \oplus \omega_r \Fa \subset\Cinf,
\] 
where $\omega = (\omega_1 : \cdots : \omega_r) \in\Omega^r\subset\BP^{r-1}(\Cinf)$ and $\Fa\subset A$ is a representative of the ideal class $[\Fa]\in\Pic(A)\cong K^\times\backslash\Af^\times/\Ahat^\times$ corresponding to the irreducible component $\Gamma_{[\Fa]}\backslash\Omega^r$ of $M_A^r(1)_{\Cinf}$ lying beneath $\Gamma\backslash\Omega^r$.

Let $t\in\GL_r(K)$ encode the given cyclic $\Fn$-isogeny, so $\Lambda_{t(\omega)}\subset\Lambda_\omega$ with
$\Lambda_\omega/\Lambda_{t(\omega)} \cong A/\Fn$. By the Invariant Factor Theorem \cite[\S22]{CurtisReiner}, there exists another basis $\omega' = (\omega'_1 : \cdots : \omega'_r)$ for $\Lambda_\omega$, such that
\[
\Lambda_{t(\omega)} = \omega'_1NA\oplus\omega'_2 A \oplus\cdots\oplus \omega'_{r-1}A\oplus\omega'_r\Fa.
\]
Now $\omega'=\gamma(\omega)$ for the base change $\gamma\in\Gamma_{[\Fa]}$ and $\gamma t\gamma^{-1}=\diag(N,1,\ldots,1)$.

After passing through the isomorphism
\[
\gamma : \Gamma\backslash\Omega^r \stackrel{\sim}{\longto} \gamma\Gamma\gamma^{-1}\backslash\Omega^r
\]
over $\Gamma_{[\Fa]}\backslash\Omega^r$, if necessary, we may assume that $t=\diag(N,1,\ldots,1)$.

The Hecke-orbit $T^\infty(x)$ contains $[H\cdot\omega]$, where $H\subset\GL_r(K)$ is the subgroup generated by $t$ and $\Gamma$.

Denote by $\delta_{ij}$ the $r\times r$ matrix with $1$ in position $(i,j)$ and $0$ elsewhere. Since $\Gamma\cap\GL_r(A)$ has finite index in $\GL_r(A)$, it follows that $G_{ij}\cap\Gamma$ has finite index in $G_{ij}$, where
\[
G_{ij} := \{1+a\delta_{ij} \;|\; a\in A\} \subset\GL_r(A),\quad G_{ij}\cong (A,+),\qquad i\neq j.
\]
Hence there exists $0\neq a\in A$ such that $1\pm a\delta_{ij}\in\Gamma$ for all $i\neq j$.

For each $n\in\BZ$ we construct the following elements of $H$:
\begin{eqnarray*}
\sigma_1^n(a) & := & t^{-n}\cdot(1+a\delta_{1r})\cdot t^n, \\
\sigma_i^n(a) & := & (1-a\delta_{i1})\cdot \sigma_1^n(-a) \cdot(1+a\delta_{i1})\cdot \sigma_1^n(a),\quad i=2,\ldots,r-1.
\end{eqnarray*}
The element $\sigma_i^n(a)\in H$ acts as translation by $a^iN^{-n}$ in the direction of the $i$th axis in $\Omega^r\subset\BA^{r-1}(\Cinf)$.

We thus see that $T^\infty(x)$ contains, in the neighborhood of every point $z\in T^\infty(x)$, points arbitrarily close to $z$ and translated in any direction. It follows that the dimension of the Zariski closure of $T^\infty(x)$ is $r-1$, and the result follows.
\end{Proof}

\begin{Cor}\label{CMdensity}
The set of CM points in $M_A^r(\CK)(\Cinf)$ is Zariski-dense. 
Likewise, if $X \subset M_A^r(\CK)_{\Cinf}$ is a special subvariety, then the CM points in $X(\Cinf)$ are Zariski-dense in $X$.
\end{Cor}

\begin{Proof}
Each irreducible component of $M_A^r(\CK)(\Cinf)$ contains at least one CM point, and its  
$T$-orbit is Zariski-dense, for a suitable Hecke correspondence $T$.
The statement for special subvarieties now follows from the definition.
\end{Proof}

\section{Modular polynomials for $A=\BF_q[T]$}\label{section4.5}

\subsection{$J$-invariants}

In this section we temporarily restrict ourselves to the case where $A=\BF_q[T]$ is a polynomial ring.

Igor Potemine \cite{Potemine} constructed a set of explicit invariants, which we denote $J_1,\ldots,J_N$, characterising the isomorphism classes of Drinfeld $\BF_q[T]$-modules, so that 
\[
M_A^r(1) = \Spec(A[J_1,\ldots,J_N]).
\] 
The value of $N$ depends on $r$ and on $q$, and the invariants are not algebraically independent.

Let $\Fn\in\BF_q[T]$ be a monic non-zero polynomial.
In \cite{BR} are constructed explicit polynomials 
\[
P_{J_i,\Fn}(X,Y_1,\ldots,Y_N) \in A[X,Y_1,\ldots,Y_N],\qquad i=1,\ldots, N,
\]
such that the graph of the Hecke correspondence $T_\Fn$ in $M_A^r(1)_K\times M_A^r(1)_K$ is defined by the equations
\[
P_{J_i,\Fn}(J'_i,J_1,\ldots,J_N) = 0,\qquad i=1,\ldots,N,
\] 
where $J_1,\ldots,J_N$ denotes the invariants from the first copy of $M_A^r(1)_K$, and $J'_1,\ldots,J'_N$ the invariants from the second copy.

It follows from \cite[Theorem 1.1]{BR} that the degree of $X$ in $P_{J_i,\Fn}(X,Y_1,\ldots,Y_N)$ is given by 
\[
\psi_r(\Fn) := |\Fn|^{r-1}\prod_{\Fp|\Fn}\frac{|\Fp|^r-1}{|\Fp|^r-|\Fp|^{r-1}}
\]
and the total degree of $P_{J_i,\Fn}(X,Y_1,\ldots,Y_N)$ is bounded by
\[
\deg P_{J_i,\Fn}(X,Y_1,\ldots,Y_N) \leq |\Fn|^{r-1}\psi_r(\Fn)^2w(J_i),
\]
where $w(J_i)>0$ is a constant depending on $J_i$. 

%

For indeterminates $X_1,\ldots,X_N$ the morphism $A[X_1,\ldots,X_N]\to A[J_1,\ldots,J_N],\; X_i\mapsto J_i$, realizes $M_A^r(1)$ as a subvariety of $\BA^N_A$.

For any subvariety $Z\subset \BA_K^N$ we define its {\em degree} $\Deg(Z)$ as the sum of the degrees of the irreducible components of the Zariski-closure of $Z$ in $\BP_K^N$ 
as projective varieties. For another subvariety $Z'\subset\BA_K^N$ we then have a form of B\'ezout's Theorem 
(\cite[Example 8.4.6]{Fulton}):
\[
\Deg(Z\cap Z') \leq \Deg(Z)\Deg(Z').
\]

Denote by $S_{\Fn}\subset \BA_{\Cinf}^{2N}$ the subvariety defined by the ideal 
\[
\langle P_{J_1,\Fn}(Y_1;X_1,\ldots,X_N), \ldots, P_{J_N,\Fn}(Y_N;X_1,\ldots,X_N)\rangle 
\subset \Cinf[X_1,\ldots,X_N,Y_1,\ldots,Y_N].
\]
By construction, the graph of the Hecke correspondence $T_\Fn$ in $M_A^r(1)_{\Cinf}\times M_A^r(1)_{\Cinf} \subset \BA_{\Cinf}^{2N}$ equals $S_{\Fn}\cap(M_A^r(1)_{\Cinf}\times M_A^r(1)_{\Cinf})$. 
Thus we have
\begin{myequation}\label{TnDeg}
\Deg(T_\Fn) = \Deg\big(S_\Fn\cap(M_A^r(1)_{\Cinf}\times M_A^r(1)_{\Cinf})\big) \leq\Deg\big(M_A^r(1)_{\Cinf}\big)^2
\prod_{i=1}^N\big(|\Fn|^{r-1}\psi_r(\Fn)^2w(J_i)\big),
\end{myequation}
where now $\Deg(T_\Fn)$ denotes the degree of the graph of $T_\Fn$ as a variety.


\subsection{Intersections}

We now drop the restriction on $A$ again.
For every transcendental element $T\in A$ we have $\BF_q[T]\subset A$ and $K/\BF_q(T)$ is an imaginary extension of degree $d_T$, say. Thus we have the canonical morphism (injective on $\Cinf$-valued points)
\[
\rho_T = f_{A/\BF_q[T]} : M_A^r(1)_K \longto M_{\BF_q[T]}^{rd_T}(1)_{K}.
\]


For any subvariety $Z\subset M_A^r(1)_K$ we define the {\em $T$-degree} $\deg_T(Z)$ to be the degree of $\rho_T(Z)\subset\BA_K^N$ as in the previous section. This leads to the following result, which we will need in the proof of Theorem \ref{MainCurves}.

\begin{Prop}\label{intersections}
Let $T\in A$ be a transcendental element, then there exist computable positive constants $c=c(T)$ and $n=n(T)$ such that the following holds. 
Let $Z\subset M_A^r(1)_{\Cinf}$ be an irreducible algebraic subvariety, and let $\FP\subset A$ be a prime which has residual degree one over $\Fp := \FP\cap\BF_q[T]$. Then
\[
\deg_T(T_\FP(Z)) \leq c\deg_T(Z)|\FP|^n.
\]
In particular, if $Z\cap T_\FP(Z)$ is finite, then
\[
|Z\cap T_\FP(Z)| \leq c\deg_T(Z)^2|\FP|^n.
\]
\end{Prop} 

\begin{Proof}
%
As $A/\FP \cong \BF_q[T]/\Fp$ we have $\rho_T(T_\FP(Z)) \subset T_\Fp(\rho_T(Z))$ ($\FP$-isogenies of Drinfeld $A$-modules are also $\Fp$-isogenies of Drinfeld $\BF_q[T]$-modules). Moreover $\dim \rho_T(T_\FP(Z)) = \dim T_\Fp(\rho_T(Z))$. Thus we have
\begin{eqnarray*}
\deg_T(T_\FP(Z)) & = & \Deg\big(\rho_T(T_\FP(Z))\big) \\
& \leq & \Deg\big(T_\Fp(\rho_T(Z))\big) \\
& = & \Deg\big(T_\Fp\cap(\rho_T(Z)\times\BA_K^N)\big) \\
& \leq & \Deg(T_\Fp)\deg_T(Z) \\
& \leq & \deg_T(Z)\Deg(M_{\BF_q[T]}^{rd_T}(1)_{\Cinf})^2
\prod_{i=1}^N\big(|\Fp|^{r-1}\psi_r(\Fp)^2w(J_i)\big),\quad\text{by (\ref{TnDeg}).}
\end{eqnarray*}
\end{Proof}

%

\section{Irreducibility of Hecke images}\label{Section5}

%
%


For a non-zero prime $\Fp\subset A$ we denote by $p_\Fp : \GL_r(\Af)\to\GL_r(K_\Fp)$ the projection onto the $\Fp$-factor. Let $\CK\subset\GL_r(\Ahat)$ be an open subgroup. Then the {\em support} of an element $g\in\GL_r(\Af)$ relative to $\CK$ is defined to be the finite set 
$\Supp_\CK(g):=\{\Fp \;|\; p_\Fp(g^{-1}\CK g) \neq p_\Fp(\CK)\}$. 
For an ideal $\Fn\subset A$ we denote by $\Supp(\Fn):=\{\Fp \;|\; \Fp|\Fn\}$ the set of prime factors of $\Fn$.

Our next goal is to prove the following result.

\begin{Thm}\label{irreducible}
Let $\CK\subset\GL_r(\Ahat)$ be an open subgroup, and let $X\subset M_A^r(\CK)_{\Cinf}$ be an irreducible Hodge generic subvariety. 
\begin{enumerate}
\item There exists a non-zero ideal $\Fm_X$, depending on $X$, with the following property. Let $g\in\GL_r(\Af)$ with $\Supp_{\CK}(g)\cap\Supp(\Fm_X)=\emptyset$ and let $\tilde{T}_g$ denote an irreducible component of $T_g$ on $M_A^r(\CK)$. Then $\tilde{T}_g(X)$ is irreducible. 
\item Suppose that $\CK\subset\GL_r(\Ahat)$ contains $d\CK := \diag\big(\det(\CK),1,\ldots,1 \big)$. 
Then there exists an open subgroup $\CK'\subset\CK$, depending on $X$, and an irreducible subvariety $X'\subset M_A^r(\CK')_{\Cinf}$ lying above $X$ such that $\tilde{T}_g(X')$ is irreducible for any irreducible component $\tilde{T}_g$ of $T_g$ on $M_A^r(\CK')_{\Cinf}$ for every $g\in\GL_r(\Af)$.
\end{enumerate}
\end{Thm}

\begin{Proof} We first replace $X$ by $X^\ns$, the non-singular locus of~$X$. 
Choose a set $S\subset\GL_r(\Af)$ of representatives for $\GL_r(K)\backslash\GL_r(\Af)/\CK$ in such a way that each $s\in S$ is of the form $s=\diag(s_0,1,\ldots,1)$ with $s_0\in\Af$. 
For each $s\in S$ let $\Gamma_s := s\CK s^{-1}\cap\GL_r(K)$.
We have $X^\an(\Cinf)\subset \Gamma_s\backslash\Omega^r$ for some $s\in S$. 
Let $\Xi\subset\Omega^r$ be an irreducible component of the preimage of $X^\an(\Cinf)$, and set $\Delta=\Stab_{\Gamma_s}(\Xi)$, so that  $X^\an(\Cinf)\cong\Delta\backslash\Xi$ (see \cite{BreuerPink}). By abuse of terminology, we call $\Delta$ the {\em analytic fundamental group} of $X$. Denote by $\hat\Delta$ and $\hat\Gamma_s$ the respective closures of $\Delta$ and $\Gamma_s$ in $\GL_r(\Af)$.

We need the following Lemma, which uses the fact that $X$ is Hodge generic.

\begin{Lem}\label{BP}
In the above situation, $\hat\Delta$ is open in $\hat\Gamma_s$.
\end{Lem}

\begin{Proofof}{Lemma \ref{BP}} 
Let $\CK'\subset\CK$ be an open subgroup for which $M_A^r(\CK')$ is a fine moduli scheme, and such that we may write
\[
M_A^r(\CK')^\an(\Cinf) = \coprod_{s\in S'}\Gamma'_s\backslash\Omega^r,
\quad\text{with $\Gamma'_s\subset\SL_r(K)$, for all $s\in S'$.}
\]
Now $X^\an(\Cinf)\cong\Delta\backslash\Xi$ has a preimage $X'^\an(\Cinf)\cong\Delta'\backslash\Xi\subset\Gamma_s'\backslash\Omega^r$ in $M_A^r(\CK')(\Cinf)$ where $\Delta',\Gamma'_s\subset\SL_r(K)$. Denote by $\hat\Delta'$ and $\hat\Gamma'_s$ the closures of $\Delta'$ and $\Gamma'_s$ in $\SL_r(\Af)$. As $X'$ is again Hodge generic, it follows from \cite[Theorem 1.1]{BreuerPink} that $\hat\Delta'$ is open in $\SL_r(\Af)$, hence $\hat\Delta'$ has finite index in $\hat\Gamma'_s$. Since $\Delta'$ has finite index in $\Delta$, and $\Gamma'_s$ has finite index in $\Gamma_s$, it follows that $\hat\Delta$ has finite index in $\hat\Gamma_s$, which proves Lemma \ref{BP}.
\end{Proofof}

We now continue with the proof of Theorem \ref{irreducible}.

{\bf(1).} 
Let $\CK_g = \CK\cap g^{-1}\CK g$.
Replacing $s$, if necessary, we may assume that the irreducible component $\tilde{T}_g$ of $T_g$ corresponds to $[s]\in\GL_r(K)\backslash\GL_r(\Af)/\CK_g$.

We have $s=\diag(s_0,1,\ldots,1)$, and we write $s_0=t_0s_1$, with $t_0\in K$ and $s_1\in\Ahat^\times$. Let $t=\diag(t_0,1,\ldots,1)\in\GL_r(K)$, then $\Delta\subset\Gamma_s\subset t\GL_r(A)t^{-1}$.

Let $\Gamma_{s,g} := s\CK_gs^{-1}\cap\GL_r(K) = sg^{-1}\CK gs^{-1}\cap\Gamma_s$, then since $\tilde{T}_g$ is given by (\ref{HeckeDiagram}), it suffices to show that $\pi^{-1}(X)$ is irreducible.
The map $\pi^{-1}(X) \to X$ is a $\Gamma_{s,g}\backslash\Gamma_s$-cover, and it suffices to show that $\Delta$ acts transitively on it. 
Since $\hat\Delta$ is open in $\hat\Gamma_s$, it follows that there exists a non-zero ideal $\Fm_X\subset A$ for which $p_\Fp(\hat\Delta) = p_{\Fp}(\hat\Gamma_s)$ for all $\Fp\nmid\Fm_X$, and consequently $\Delta$ and $\Gamma_s$ have the same image in $t\GL_r(A/\Fn)t^{-1}$ for all ideals $\Fn\subset A$ prime to $\Fm_X$. 

Let $\Fn = \prod_{\Fp\in\Supp_\CK(g)}\Fp^{n_\Fp}$, with each $n_{\Fp}$ sufficiently large that 
\begin{myequation}\label{smallkernel}
\ker\big(p_\Fp(s\CK s^{-1})\to t\GL_r(A/\Fp^{n_\Fp})t^{-1}\big) \subset p_\Fp(sg^{-1}\CK gs^{-1}).
\end{myequation}
(Clearly (\ref{smallkernel}) holds with $n_\Fp=0$ when $\Fp\not\in\Supp_\CK(g)$.)
Then $\Fn$ is prime to $\Fm_X$, and we see that
\[
\ker\big(\Gamma_s\to t\GL_r(A/\Fn)t^{-1}\big) \subset \Gamma_{s,g}.
\] 
It follows that $\Gamma_{s,g}\backslash\Gamma_s \hookrightarrow t\GL_r(A/\Fn)t^{-1}$, with the same image as $\Delta$, and hence $\Delta$ acts transitively on $\Gamma_{s,g}\backslash\Gamma_s$, as desired.

\medskip

{\bf (2).} 
We continue using the above notation. 

Replacing $\CK$ by its intersection with $\{x\in\GL_r(\Ahat) \;|\; \det(x)\equiv 1 \bmod \Fn\}$ for a suitable $\Fn\subset A$, if necessary, we may assume that $\Gamma_s\subset\SL_r(K)$, while $d\CK\subset\CK$ still holds.

We set
\[
\CK' := s^{-1}\hat\Delta \cdot d\CK s = s^{-1}\hat\Delta s \cdot d\CK, 
\] 
and we must first show that $\CK'$ is an open subgroup of $\CK$. Since $s^{-1}\hat\Delta s \subset \CK$, we have $\CK'\subset\CK$, and it remains to show the index is finite. We embed $\Ahat^\times\into\GL_r(\Ahat)$ via $x\mapsto\diag(x,1,\ldots,1)$, so $\GL_r(\Ahat) = \SL_r(\Ahat)\cdot \Ahat^\times$. Next, $s\hat\Delta s^{-1}$ has finite index in $\SL_r(\Ahat)$ by \cite[Theorem 1.1]{BreuerPink}, and $\det(\CK)$ has finite index in $\Ahat^\times$, so it follows that $\CK'$ has finite index in $\GL_r(\Ahat)$.


Now let $\Gamma_s' := \Gamma_s\cap s\CK's^{-1} = s\CK's^{-1}\cap\SL_r(K) = \hat\Delta\cap\SL_r(K)$.
Choose an irreducible $X'\subset M_A^r(\CK')_{\Cinf}$ lying above $X$ and below $\Xi$. Its analytic fundamental group is $\Stab_{\Gamma'_s}(\Xi) = \Delta\cap\Gamma_s' = \Delta\cap\hat\Delta\cap\SL_r(K) = \Delta$. 

Since $\hat\Delta$ is open in $\SL_r(\Af)$ and $\SL_r(K)$ is dense in $\SL_r(\Af)$, it follows that the closure of $\Gamma_s'$ in $\SL_r(\Af)$ equals $\hat\Delta$. The argument from part (1) now shows that $\tilde{T}_g(X')$ is irreducible, and there are no ``bad'' primes that need to be excluded from $\Supp_\CK(g)$. The result follows.
\end{Proof}

\section{Proof of the main results}\label{MainProofsSect}

\subsection{CM points with good behaviour at a given prime}
%
%
%
%

\begin{Proofof}{Theorem \ref{MainMoonen}}
Suppose $X\subset M_A^r(\CK)_{\Cinf}$ is an irreducible algebraic subvariety, and $X(\Cinf)$ contains a Zariski-dense subset $\Sigma$ of CM points satisfying conditions (1) and (2) of Theorem \ref{MainMoonen}. Clearly it suffices to prove the Theorem for $\CK=\GL_r(\Ahat)$. 

For any CM point $x$ we will write $R_x:=\End(\phi^x)$, which has conductor $\Fc_x$.

Let $(A',h)$ be a Hodge pair for $X$, so $X$ is an irreducible subvariety of $T_h(M_{A'}^{r'}(1)_{\Cinf})$ in $M_A^r(1)_{\Cinf}$. Let $X'\subset M_{A'}^{r'}(1)_{\Cinf}\subset M_A^r(1)_{\Cinf}$ be an irreducible subvariety such that $X$ is an irreducible component of $T_h(X')$. $X'$ is Hodge generic in $M_{A'}^{r'}(1)_{\Cinf}$. Denote by $\Sigma'\subset X'(\Cinf)$ a set of CM points in $X'(\Cinf)$ such that $\Sigma\subset T_h(\Sigma')$; it is Zariski-dense in $X'$. Furthermore, the $\Fp$-order of $\Fc_x$ is bounded independently of $x\in\Sigma'$.
Let $\Fp'\subset A'$ be a prime above $\Fp$ which is residual in every $\Quot(R_x)$, it exists as $A'\subset\cap_{x\in\Sigma'}\Quot(R_x)$, though we may have to replace $\Sigma'$ by a Zariski-dense subset. Replace $(X,\Sigma,A,r,\Fp)$ by $(X',\Sigma',A',r',\Fp')$, and we have reduced ourselves to the case where $X$ is Hodge generic.

The $\Fp$-orders of the conductors $\Fc_x$ are bounded by $m$. It follows that, for every $x\in\Sigma$, there exists an $A$-submodule $H_x\subset (A/\Fp^m)^r$ and an isogeny $f_x : \phi^x\to\psi^x$ with $\ker(f_x)\cong H_x$ such that $\Fp$ does not divide the conductor of $\End(\psi^x)$. Since $(A/\Fp^m)^r$ has only finitely many $A$-submodules, we may, after replacing $\Sigma$ by a Zariski-dense subset, assume that all these modules $H_x$ are equal to a fixed $H\subset (A/\Fp^m)^r$. 
Let $h\in\GL_r(\Af)$ such that 
$T_h$ represents isogenies with kernels isomorphic to $H$, and so for every $x\in\Sigma$, the Hecke image $T_h(x)$ contains at least one CM point $y$ such that $\Fp\nmid\Fc_y$. A subset $\Sigma'$ of these $y$'s is again Zariski-dense in an irreducible component $X'$ of $T_h(X)$. It suffices to show that this $X'$ is special. Since $X'$ is again Hodge generic, we may replace $(X,\Sigma)$ by $(X',\Sigma')$, and have now reduced to the situation where $X$ is Hodge generic and contains a Zariski-dense set $\Sigma$ of CM points $x$ for which $\Fp$ is residual in $R_x$. 

Let $\CK'\subset\GL_r(\Ahat)$ and $X'\subset M_A^r(\CK')_{\Cinf}$ be the level structure and preimage of $X$ given by Theorem \ref{irreducible}.2. The set $\Sigma$ lifts to a Zariski-dense subset $\Sigma'\subset X'(\Cinf)$ of CM points with the same endomorphism rings. It suffices to show that $X'$ is an irreducible component of $M_A^r(\CK')_{\Cinf}$, so, again, we replace $(X,\Sigma)$ by $(X',\Sigma')$.

\noindent\begin{minipage}{12cm}
\hspace{0.6cm}Let $F$ be a field of definition of $X$. 
Since $\Sigma\subset X(\bar{K})$, we may choose $F$ to be a finite extension of $K$. Denote by $F_s$ the separable closure of $K$ in $F$. Let $n\in\BN$ be sufficiently large that $\Fp^n$ is principal $K$, and $\FP^n$ is principal for every prime $\FP$ above $\Fp$ in the Hilbert class field of $F$. 

\hspace{0.6cm}Let $x\in\Sigma$, $K'=\Quot(R_x)$ and $K'(x)$ the field of definition of $x$ over $K'$ from Proposition \ref{CMmain}.
Let $L$ be the Galois closure of $F_sK'(x)$ over $K'$, and let $\FP\subset K'$ be a prime above $\Fp$ with residue degree $f(\FP|\Fp)=1$.
Denote by $\sigma\in\Aut(FL/FK')$ an extension of the inverse Artin element $(\FP^n,L/K')^{-1}$, it fixes $F$ since $\Fp^n$ is principal in the Hilbert class field of $F$. 

\hspace{0.6cm}Now by Proposition \ref{CMmain} we have $\sigma(x)\in T_{\Fp^n}(x)$. 
Since $\sigma$ fixes $F$, a field of definition for $X$ and $T_{\Fn}(X)$, we have
%
%
%
%
\end{minipage}
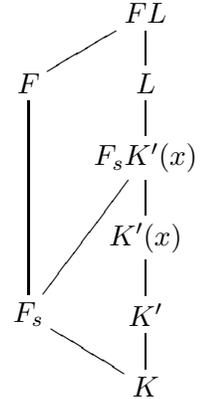
\begin{minipage}{4cm}
\[
\xymatrix@=13pt{
& FL\ar@{-}[d]\ar@{-}[dl] \\
F\ar@{-}[ddd] & L\ar@{-}[d] \\
& F_sK'(x)\ar@{-}[d]\ar@{-}[ddl] \\
& K'(x)\ar@{-}[d] \\
F_s\ar@{-}[dr] & K'\ar@{-}[d] \\
& K
}
\]
\end{minipage}
\begin{myequation}\label{lotsofx}
x\in X(\Cinf)\cap T_{\Fp^n}(X(\Cinf)).
\end{myequation}

%

Restricting $\Sigma$ to a Zariski-dense subset, if necessary, it follows from (\ref{lotsofx}) that for a fixed irreducible component $\tilde{T}_{\Fp^n}$ of $T_{\Fp^n}$, we have
\[
\Sigma\subset X(\Cinf)\cap \tilde{T}_{\Fp^n}(X(\Cinf)).
\]
But $\tilde{T}_{\Fp^n}(X)$ is irreducible by Theorem \ref{irreducible}.2, so $X=\tilde{T}_{\Fp^n}(X)$.

Notice that $\tilde{T}_{\Fp^n}$ preserves the irreducible component of $M_A^r(\CK)_{\Cinf}$ containing $X$. 
Now by Theorem \ref{density}, the $\tilde{T}_{\Fp^n}$-orbit of any point in $X(\Cinf)$ is Zariski-dense in an irreducible component of $M_A^r(\CK')_{\Cinf}$, but this orbit lies completely inside $X(\Cinf)$, and the result follows.

Finally, suppose that $\Sigma$ lies in one Hecke orbit. Then the $K_x$ are all equal to a fixed field $K'$, and the set of prime factors of the conductors $\Fc_x$ is finite. Denote by $K'_s$ the separable closure of $K$ in $K'$. Then, by the \v{C}ebotarev Density Theorem \cite[Theorem 5.6]{Fried-Jarden}, there exist infinitely many primes $\Fp\subset A$ which split completely in (the Galois closure of) $K'_s$, hence are residual in $K'$, and do not divide any $\Fc_x$. Conditions (1) and (2) are thus satisfied and $X$ is special.
\end{Proofof}

\subsection{CM points on curves}
%


%

We now prove our second main result.

\medskip

%

\begin{Proofof}{Theorem \ref{MainCurves}}
It suffices to prove the Theorem for $\CK = \GL_r(\Ahat)$. Let $(A',h)$ be a Hodge pair for $X$, then we may replace $X$ by a preimage in $M_{A'}^{r'}(1)_{\Cinf}$, so we may assume that $X$ is Hodge generic. Our goal is to show that $X$ is an irreducible component of $M_A^r(1)_{\Cinf}$.

Let $F$ be a field of definition for $X$ containing the Hilbert class field of $K$, we may assume that $F$ is a finite extension of $K$, as $X(\bar{K})$ is infinite. Let $\Fm_X \subset A$ be the ideal given in Theorem \ref{irreducible}. 
Fix a transcendental element $T\in A$ such that $K/\BF_q(T)$ is a separable geometric extension, 
and let $c=c(T)$ and $n=n(T)$ denote the constants from Proposition \ref{intersections}.

Let $x\in X(\bar{K})$ be a CM point. Suppose that there exists a non-zero principal prime $\Fp\subset A$ with the following properties:
\begin{enumerate}
\item[(1)] $\Fp\nmid\Fm_X$.
\item[(2)] $\Fp$ is residual in $R=\End(\phi^x)$ and in $F$.
\item[(3)] $\Fp$ has residual degree one over $P:=\Fp\cap\BF_q[T]$.
\item[(4)] $|\Pic(R)| > c[F:K]\deg_T(X)^2|\Fp|^n$.
\end{enumerate}

We want to show that $X\subset T_\Fp(X)$. Suppose that the intersection $X\cap T_\Fp(X)$ is proper, then $|X\cap T_\Fp(X)| \leq c\deg_T(X)^2|\Fp|^n$, by Proposition \ref{intersections}. 

Denote by $F_s$ the separable closure of $K$ in $F$, let $K'=\Quot(R)$ and let $K'(x)$ be the ring class field of $R$, which is a field of definition of the point~$x$. Let $L$ be the Galois closure of $F_sK'(x)$ over $K'$. Let $\FP$ be a prime of $L$ above $\Fp$ and let $\sigma\in\Aut(FL/FK')$ be an extension of the Artin element $(\FP, L/K')^{-1}$ (it also fixes $F$, as $\Fp$ is residual in $F$). Then $\sigma(x)\in T_\Fp(x)$ by Proposition \ref{CMmain}, and $\sigma$ fixes $F$, which is a field of definition of $X$ and $T_\Fp(X)$. Therefore, we have 
\[
\sigma(x)\in X^{\sigma}(\Cinf)\cap T_\Fp(X(\Cinf)) = X(\Cinf)\cap T_\Fp(X(\Cinf)),
\]
and furthermore
\[ 
\Gal(FK'(x)/FK')\cdot x\subset X(\Cinf)\cap T_\Fp(X(\Cinf)). 
\] 
But
\[
|\Gal(FK'(x)/FK')\cdot x| \geq |\Pic(R)|/[F:K],
\]
and by assumption (4), this is larger than the intersection degree of $X\cap T_\Fp(X)$. Therefore this intersection is improper and we have $X\subset T_\Fp(X)$. 
Since $T_\Fp(X)$ is irreducible by Theorem \ref{irreducible}.1, it follows that $X=T_\Fp(X)$, and so $X(\Cinf)$ contains the entire $T_\Fp$-orbit of any point. But such orbits are dense in an irreducible component of $M_A^r(1)_{\Cinf}$, by Theorem \ref{density}, so we are done.

\medskip

Since there are only finitely many CM points $x$ with $\HCM(x)$ smaller than a given constant (Proposition \ref{CMheights}.1), it remains to show that, if $\HCM(x)$ is sufficiently large, then a prime
$\Fp$ satisfying (1)-(4) above does indeed exist. 

Let $K'_s$ denote the separable closure of $\BF_q(T)$ in $K'$, let $H_A$ denote the ring class field of $A$ and denote by $M$ the Galois closure of $F_sK'_sH_A$ over $\BF_q(T)$. For $t\in\BN$ we denote
\[
\pi_M(t) := |\{P \in\BF_q[T] \;|\; \text{$P$ is prime, $P$ splits in $M$ and $|P|=q^t$}\}|.
\]
Denote by $n_c$ the constant extension degree of $M/\BF_q(T)$ and by $n_g$ the geometric
extension degree. We have $n_gn_c = [M:\BF_q(T)] \leq C_1$ for some constant
$C_1$ not depending on $x$. Now the strong \v{C}ebotarev Theorem for function fields 
\cite[Prop. 5.16]{Fried-Jarden} says:
\[
\text{If $n_c|t$, then} \quad |\pi_M(t) - \frac{1}{n_g}q^t/t| < 4(g_M + 2)q^{t/2},
\]
where $g_M$ denotes the genus of $M$. (Notice that this strong version of the \v{C}ebotarev Theorem is where one requires the Generalized Riemann Hypothesis in characteristic 0). We can
bound $g_M$ in terms of the genus $g'$ of $K'$ (which equals the genus of $K_s$) using the Castelnuovo inequality
\cite[III.10.3]{Stichtenoth}, and obtain
\[
\pi_M(t) > \frac{1}{C_1}q^t/t - (C_2g'+C_3)q^{t/2},
\]
where $C_1,C_2$ and $C_3$ are positive constants independent of our CM point 
$x$ and of $t$. 
Let $\Fc$ be the conductor of $R$ in $K'$. A principal prime $\Fp$ satisfying (1), (2)
and (3) above will exist if $\pi_M(t) > \log|\Fm_X| + \log|\Fc|$, where 
$|\Fp| = |P| = q^t$, because $P$ splits completely in $F_sK'_sH_A/\BF_q(T)$ and ramifies totally in $FK'/F_sK'_s$, hence has residue degree 1 in $FK'/K$.

On the other hand, we also want to satisfy (4). For this we will use Proposition \ref{CMheights}.2:
\begin{eqnarray*}
|\Pic(R)| & > & C_\varepsilon\HCM(R)^{1-\varepsilon} 
\quad\forall\varepsilon >0\\
& = & C_\varepsilon\big(q^{g'}|\Fc|^{1/r}\big)^{1-\varepsilon},
\end{eqnarray*}
which in turn needs to be larger than $c[F:K]\deg_T(X)^2|\Fp|^n$. 

In summary, a suitable prime $\Fp \subset A$ will exist if there is a
simultaneous solution $t\in n_c\BN$ to the two inequalities:
\[
\frac{1}{C_1}q^t/t - (C_2g'+C_3)q^{t/2} > \log|\Fm_X|+\log|\Fc|,
\]
%
\[
C_{\varepsilon}\big(q^{g'}|\Fc|^{1/r}\big)^{1-\varepsilon} > 
c[F:K]\deg_T(X)^2q^{tn},\quad\text{for some $\varepsilon >0.$}
\]

(Note that $n_c$ does depend on $x$, but it is bounded by $C_1$).
If $\HCM(x) = q^{g'}|\Fc|^{1/r}$ is sufficiently large, then such a solution 
does exist, as one just needs $q^t$ to be large compared to $g'$ and $\log|\Fc|$, and small compared to $q^{g'}$ and $|\Fc|$. This completes the proof.
\end{Proofof}

\begin{Rem} 
{\em The proof of Theorem \ref{MainCurves} is actually effective. Suppose that
$X$ is defined over a finite extension $F/K$. Then in order to prove
that a curve $X$ is special, it suffices to find {\em one} CM point $x$ with
$\HCM(x)$ larger than some computable constant depending on $F$, $T$, $\deg_T(X)$, and on $\Fm_X$. 
}
\end{Rem}

{\small


\bigskip

\parbox{10cm}{
Department of Mathematical Sciences\\
University of Stellenbosch\\
Stellenbosch 7602\\
South Africa\\
fbreuer@sun.ac.za\\
tel: +27 21 8083288\\
fax: +27 21 8083828}

}

\end{document}